\documentclass{article}
\usepackage{amsmath}
\usepackage{amsfonts}
\usepackage{tikz,tkz-tab}
\usetikzlibrary{matrix,arrows,decorations.pathmorphing} 
\newtheorem{lemma}{Lemma}
\newtheorem{theorem}{Theorem}

\newtheorem{claim}{Claim}
\begin{document}
\title{The number of binary rotation words}

\author{A. Frid\thanks{Sobolev Institute of Mathematics SB RAS and Universit\'e de Lorraine; supported in part by Presidential grant MK-4075.2012.1 and by RFBR grant 12-01-00089.}, 
D. Jamet\thanks{Loria and Universit\'e de Lorraine.}}

\maketitle
%
\begin{abstract}
We consider binary rotation words generated by partitions of the unit circle to two intervals and give a precise formula for the number of such words of length $n$. We also give the precise asymptotics for it, which happens to be $\Theta(n^4)$. The result continues the line initiated by the formula for the number of all Sturmian words obtained by Lipatov in 1982, then independently by Berenstein, Kanal, Lavine and Olson in 1987, Mignosi in 1991, and then with another technique by Berstel and Pocchiola in 1993.
\end{abstract}
%
%
%

-----
-----






\section{Introduction}
Infinite words arising from rotations of the circle belong to the same family of infinite words defined by the means of dynamical systems as Sturmian words and interval exchange words in general. They were considered by G. Rote in 1992 \cite{rote} and can be defined using three parameters $\alpha,\beta,\gamma \in [0,1)$ as $r=r_0r_1\cdots$, where for all $i$ we have
\begin{equation}\label{e:def}
r_i=\begin{cases} 1, \mbox{~if~} \{i\alpha\} \in [\beta,\gamma), \\ 0, \mbox{~otherwise}. \end{cases}
\end{equation}
(Here the interval $[\beta,\gamma)$ is denoted as usual if $\beta< \gamma$ and as $[\beta;1)\cup [0,\gamma)$ otherwise.)

In the particular case when $\gamma-\beta = \alpha \pmod{1}$, $w$ is a {\it Sturmian} word. The family of Sturmian words is very well studied (see Chapter 2 of \cite{lothaire_ch2}); in particular, the total number of factors of all Sturmian words taken together is known to be
\[1+\sum_{p=1}^{n}(n-p+1)\varphi(p),\]
where $\varphi$ is the Euler's totient function. This formula was rediscovered several times \cite{lipatov,bklo,mignosi,bp1}; the order of growth of this function is $\Theta(n^3/\pi^2)$.

In \cite{cf,fr_st2} Cassaigne and the first author estimated and for some cases found the number of factors of length $n$ of all rotation words with a given length $\gamma-\beta$ of the interval; it happens that it also grows as $\Theta(n^3)$. In \cite{afmp}, Ambro\v{z}, Mas\'akov\'a, Pelantov\'a and the first author estimated the number of all words arising from three-interval exchange, which continues the same line since Sturmian words are exactly two-interval exchange words; it happens that the number of three-interval exchange words grows as $\Theta(n^4)$. In \cite{bv}, Berstel and Vuillon coded rotation words by Sturmian words.

In this paper, we find a precise formula for the number of all rotation words \eqref{e:def}, predictably involving sums of the Euler's function. To write down the formula, we had to understand very clearly the structure of the set of rotation words, which is of independent interest.

\section{Main statement}
The main result of the paper is the following
\begin{theorem}\label{t:main}
Starting from $n=3$, the number of binary rotation words of length $n+1$ is
\begin{equation}\label{e:main}
f(n+1)=n^2+3n+4+\frac{1}{2}\sum_{p=3}^n \varphi(p)(n^2-p^2+n+p)-f_1(n)-2\sum_{l=2}^{n-1} f_2(n,l),\end{equation}
where 
\begin{equation}\label{e:f1}
f_1(n)=\begin{cases} 2 \sum_{i=k}^{2k}\sum_{p=1}^{i+1}\varphi(p), \mbox{~if~} n=2k+1, \\  
2 \sum_{i=k}^{2k-1}\sum_{p=1}^{i+1}\varphi(p) + \sum_{p=1}^{k}\varphi(p), \mbox{~if~} n=2k, \end{cases}
\end{equation}
\begin{equation}\label{e:g}
g(n,l)= n-l+1+(n \bmod (l+1)),
\end{equation}
\[h(n,l)=\min(l+1,n-l),\]
and 
\begin{equation}\label{e:f2}
 f_2(n,l)=\left ( \frac{1}{2}\left \lfloor \frac{n}{l+1} \right\rfloor g(n,l)-h(n,l)\right )(\varphi(l+1)-1)+h(n,l)
\left ( \frac{\varphi(l+1)}{2}-1 \right ).
\end{equation}
\end{theorem}
Note that the only addend of this formula growing faster than than $O(n^3)$ is the sum 
\[\sum_{p=3}^n \varphi(p)(n^2-p^2).\]
So, the asymptotics of the number of binary rotation words is equal to the asymptotics of this addend, which means that
\[f(n)= \frac{3 n^4}{4 \pi^2}+O(n^3\log n).\]
The values of $f(n)$ for some values of $n$ are shown in the table below. 

{\small
\begin{center}
\begin{tabular}{c|c|c|c|c|c|c|c|c|c|c}
 $n$ & 6&7&10&15&20&30&50&75&100 \\ \hline
$f(n)$ &  64&112&504&2804&9442&51306&423814&2222984&7155096 \\ \hline
\rule{0pt}{1.62em}$\displaystyle \frac{4\pi^2f(n)}{3 n^4} \approx$
&0.65&0.61&0.66&0.73&0.78&0.83&0.89&0.92&0.94
\end{tabular}
\end{center}
}

\medskip
The rest of the paper is devoted to the proof of Theorem \ref{t:main} and thus to a study of internal structure of the set of rotation words.

\section{Rotations and Sturmian words}
Denote the prefix $r_0\cdots r_{n-1}$ of length $n$ of the word $r$ defined in \eqref{e:def} by $r(\alpha,\beta,\gamma,n)$. The parameter $\alpha$ is called the {\it slope} of the rotation word $r$. The set of all rotation words $r(\alpha,\beta,\gamma,n)$ of length $n$ is denoted by $R(n)$, so, the searched function is $f(n)=\#R(n)$. 
\begin{lemma}
It is sufficient to consider rotation words of slopes not greater than 1/2:
\[R(n)=\{r(\alpha,\beta,\gamma,n)|\alpha \in (0,1/2), \beta,\gamma \in \mathbb{R}/\mathbb{Z}.\}\]
\end{lemma}
{\sc Proof.} Due to the symmetry, we have $r(\alpha,\beta,\gamma,n)=r(1-\alpha,1-\beta,1-\gamma,n)$ if $\{k\alpha\}\neq \beta$ or $\gamma$ for all $k=0,\ldots,n-1$, that is, if the point $k\alpha$ is never equal to the end of an interval. But if it is, we can just take $r(1-\alpha,1-\beta,1-\gamma,n)$ and then slightly shift the interval to avoid its ends. So, slopes less than 1/2 and greater than 1/2 give exactly the same set of all rotation words. \hfill $\Box$

\medskip
The following lemma is a particular case of the result of Berstel and Vuillon \cite{bv}. We give its proof for the sake of clarity.
\begin{lemma}
For any two-interval rotation word $r=r(\alpha,\beta,\gamma,n)$, where $\alpha \leq 1/2$, we have
\begin{equation}\label{rk}
r_k=r_{k-1}+u_k-v_k
\end{equation}
for the Sturmian words $u=r(\alpha,\beta,\beta+\alpha,n)$ and $v=r(\alpha,\gamma,\gamma+\alpha,n)$ of the slope $\alpha$.
\end{lemma}
{\sc Proof.} The fact that $u_k=1$ is equivalent to the fact that $\beta \in [\{(k-1)\alpha\},\{k\alpha\})$; the fact that $v_k=1$ is equivalent to the fact that $\gamma \in [\{(k-1)\alpha\},\{k\alpha\})$. So, if $u_k=v_k=0$, the interval $[\{(k-1)\alpha\},\{k\alpha\})$ contains neither $\beta$ nor $\gamma$, and thus $r_k=r_{k-1}$; if $u_k=v_k=1$, the interval $[\{(k-1)\alpha\},\{k\alpha\})$ contains both $\beta$ and $\gamma$, and thus $r_k=r_{k-1}$ again; if $u_k=1$ and $v_k=0$, then $[\{(k-1)\alpha\},\{k\alpha\})$ contains $\beta$ but not $\gamma$, and thus $r_{k-1}=0$ and $r_k=1$; at last, if $u_k=0$ and $v_k=1$, then $[\{(k-1)\alpha\},\{k\alpha\})$ contains $\gamma$ but not $\beta$, and thus $r_{k-1}=1$ and $r_k=0$. In all the four cases \eqref{rk} holds.  \hfill $\Box$

\medskip
Note that the symbols $u_0$ and $v_0$ are not used in the previous lemma, so, we see that a rotation word of length $n+1$ is uniquely defined by its first symbol and two Sturmian words of the same slope of length $n$. If these two Sturmian words are distinct, we can uniquely reconstruct from them the symbol $r_0$; if they are equal, both rotation words $0^{n+1}$ and $1^{n+1}$ can appear. It is clear also that each pair of Sturmian words of the same slope gives some rotation word (of that slope). This gives us the next lemma:

\begin{lemma}
The number of binary rotation words is bounded as
\begin{equation}\label{pair_st} 
f(n+1)\leq \#\{(u,v)|u,v \in St(n,\alpha), \alpha \in (0,1/2), u\neq v\}+2.~~~~~ \Box
\end{equation}
\end{lemma}
Here $St(n,\alpha)$ is the set of all Sturmian words of length $n$ and of slope $\alpha$. The addend 2 in the formula above corresponds to all possible pairs of equal Sturmian words which all correspond to two rotation words, $0^{n+1}$ and $1^{n+1}$.

Denote by $f_{pairs}(n)$ the number of such pairs of distinct Sturmian words of length $n$ of the same slope $\alpha \in (0,1/2)$, so that \eqref{pair_st} can be rewritten as
\begin{equation}\label{e:upperbound1}
f(n+1)\leq f_{pairs}(n)+2.
\end{equation}

\begin{lemma} \label{l:pairs} 
For all $n\geq 1$ we have
\[f_{pairs}(n)=n(n+1)+\frac{1}{2}\sum_{p=3}^n \varphi(p)(n^2-p^2+n+p).\]
\end{lemma}
{\sc Proof.} Recall that the Farey series $F_n$ of order $n$ is the increasing sequence of all fractions between 0 and 1 whose denominators are at most $n$. The intervals between consecutive Farey fractions are called {\it Farey intervals}. The first Farey fraction is taken to be $0=0/1$; all the others are of the form $q/p$, where $1\leq q < p\leq n$, gcd$(q,p)=1$.

It is well-known that the sets $St(n,\alpha)$ coincide for all $\alpha$ from the same Farey interval of order $n$; if the beginning of the interval is the fraction $q/p$, we can denote this set as $St(n,\alpha)=St(n,q,p)$.

Let us say that a Sturmian word of length $n$ is {\it new} in the Farey interval starting from $q/p$ if it belongs to $St(n,q,p)$ but does not belong to any $St(n,q',p')$ for $q'/p' < q/p$, where $q/p, q'/p' \in F_n$. Denote the set of all new Sturmian words from $St(n,q,p)$ by $N(n,q,p)$; all other words from $St(n,q,p)$ are called {\it old}, and their set is denoted by $Old(n,q,p)$.

As it follows directly from the results by Berstel and Pocchiola \cite{bp2},
for all $q/p \in F_n$ with $p>1$ we have
 $\#N(n,q,p)=n-p+1$.
As a corollary, we immediately see that $\#Old(n,q,p)=p$.

Now let us count $f_{pairs}(n)$ starting from the minimal slope and going on along the Farey series. In the interval starting from 0, all the words are new, and they give $n(n+1)$ pairs. In any other interval, we are interested only in pairs where at least one of the words is new, since the pairs where both words are old have been counted before. So, after excluding pairs of old words, we see in the interval starting from $q/p$ the following number of new Sturmian pairs:
\[n(n+1)-\#Old(n,q,p)(\#Old(n,q,p)-1)=n^2-p^2+n+p.\]
Now note that the number of Farey fractions whose denominator is $p$ and which are less than 1/2 is $\varphi(p)/2$ for all $p>2$; for $p=2$, the only Farey fraction is 1/2, but we are not interested in slopes greater than 1/2; for $p=1$, the case of $0=0/1$ is a bit special and has been considered in the beginning of this paragraph. So, summing up, we obtain that
\[f_{pairs}(n)=n(n+1)+\frac{1}{2}\sum_{p=3}^n \varphi(p)(n^2-p^2+n+p). ~~~~~\Box\]

Together with \eqref{e:upperbound1}, this lemma already gives us an upper bound for $f(n+1)$.
However, to pass to a precise formula, we should classify the cases when different pairs of Sturmian words give the same rotation word.

We start from the following
\begin{lemma}\label{0011}
If a rotation word $r$ contains both factors 00 and 11, then it appears from only one pair of Sturmian words $u$ and $v$ of the same slope $\alpha<1/2$. They can be found by the equalities $u_k=v_k=0$ if $r_{k-1}=r_k=1$ and $u_k=r_k$, $v_k=r_{k-1}$ otherwise.
\end{lemma}
{\sc Proof.} Let $a$ be the symbol whose interval is shorter than the other one and thus not longer than 1/2. The jump $\alpha$ of the moving point cannot be greater than the length of the longer interval, and thus, since $aa$ appears in $r$, we see that in any pair of Sturmian words of the slope $\alpha<1/2$, generating $r$, the slope $\alpha$ is less than the length of the interval corresponding to $a$, and all the more of the other interval. It means exactly that we can never have $u_k=v_k=1$. So, we can uniquely reconstruct the words $u$ and $v$: if $r_k=r_{k-1}$, then $u_k=v_k=0$, if $r_k=1$ and $r_{k-1}=0$, then $u_k=1$ and $v_k=0$, and if $r_k=0$ and $r_{k-1}=1$, then $u_k=0$ and $v_k=1$. This is equivalent to the statement of the lemma.  \hfill $\Box$

\medskip
So, to classify all pairs of Sturmian words of the same slope $\alpha$, we must consider only those of them which contain consecutive occurrences of at most one symbol. Due to the symmetry, we can suppose for a while that this symbol is 0, that is, that the rotation words considered do not contain the factor $11$.

The words $0^{n+1}$ and $1^{n+1}$ have been already excluded from consideration and gave the addend 2 to the formula \ref{e:upperbound1}. So, in what follows we consider two cases: either $r$ contains only one symbol 1, that is, $r=0^i10^{n-i}$ for some $i$, or $r$ contains at least two 0s and two 1s. The proof will rely on properties of Sturmian words, so that we start from listing some of them.

\section{Properties of Sturmian words}\label{s:prop}
Denote the set of all Sturmian words of the slope $\alpha$ by $St(\alpha)$; so, $St(\alpha,n)=\{0,1\}^n \cap St(\alpha)$. As always, we say that a word $u$ from a binary language $L$ is left (right) {\em special} in $L$ if both $0u$ and $1u$ (respectively, $u0$ and $u1$) are also in $L$. The mirror image $\tilde{u}$ of the word $u=u_1\cdots u_m$ is the word $u_m \cdots u_1$.

The following statements are classical and can be found in any survey on Sturmian words (see, e.~g., \cite{lothaire_ch2}).
\begin{claim}\label{cl:st_sp}
Each language $St(\alpha)$ contains exactly one left special word $u$ and one right special word $v$ of each length, and $u=\tilde{v}$. A shorter left special word from $St(\alpha)$ is always a prefix of a longer one; the same holds symmetrically on right special words and suffixes.
\end{claim}
\begin{claim}\label{cl:st_numbersp}
 The total number of left (right) special Sturmian words of length $n$ is $\sum_{p=1}^{n+1}\varphi(p)$.
\end{claim}

Another family of known facts concerns the construction of Sturmian words with directive sequences, standard words and central words. The facts below can be found in \cite{lothaire_ch2}.

Consider a {\it directive sequence} $(d_1,d_2,\cdots)$, where $d_1\geq 0$, $d_i>0$ for all $i>1$, and construct a sequence of words
\[s_{-1}=1, s_0=0, s_n=s_{n-1}^{d_n} s_{n-2} (n \geq 1).\]

The words $s_n$ obtained are called {\it standard} words. All standard words are Sturmian, and for each Sturmian word $w$ there exists a directive sequence such that all the standard words it generates are factors of $w$. So, the set of factors of $w$ is uniquely determined by some directive sequence, which is directly related to the continuous fraction expansion $[0,1+d_1,d_2,\cdots]$ of the slope $\alpha$ of $w$. In what follows we denote this slope by $\alpha_d$. 

The slope $\alpha_d$ is not greater than 1/2 if and only if $d_1 >0$. Since in this paper we are interested in Sturmian words whose slope is less than 1/2, from now on we assume for each directive sequence that $d_i>0$ for all $i >0$.

It can be easily checked that all standard words of length at least 2 end by 01 or 10. For each standard word $s_n=c_n ab$, where $a\neq b$, $a,b \in \{0,1\}$, the word $c_n$ is called a {\it central} word. In what follows we denote the word $c_n ba$ by $s_n'$.

The following facts on standard and central words can be easily proved. Once again, we refer to \cite{lothaire_ch2} for details.

\begin{claim}
All standard words from $St(\alpha_d)$ are left special in that language.
\end{claim}
\begin{claim}
For each standard word $s_n \in St(\alpha_d)$, the word $s_n'$ is also standard and belongs to $St(\alpha_d)$, but is not left special in that language. 
\end{claim}
\begin{claim}\label{c:0110}
 If $s_n=s_{n-1}^{d_n} s_{n-2}$, then $s_n'=s_{n-1}^{d_n-1} s_{n-2} s_{n-1}$.
\end{claim}
\begin{claim}
Central words $c_n$ are bispecial in $St(\alpha_d)$.
\end{claim}
\begin{claim}
A word $c$ is bispecial in $St(\alpha_d)$ if and only if it is obtained by deleting two last symbols from some word $s_{n-1}^{t}s_{n-2}$, where $0<t\leq d_n$. We shall denote this word by $c_{n,t}$; in particular, $c_{n,d_n}=c_n$.
\end{claim}
\begin{claim}
 The total number of central words of a length $n$ is $\varphi(n+2)$.
\end{claim}

\begin{claim}
The last two symbols of standard words alternate: if $s_{n-1}$ ends by $01$, then $s_n$ ends by $10$, and vice versa.
\end{claim}

For each directive sequence $(d_1,d_2,\ldots)$, we denote the length of the $n$th standard word $s_n$ by $l_n$.

The following lemma can be proved analogously to Theorem 2.2.31 from \cite{lothaire_ch2}.
\begin{lemma}\label{dn+2}
 For all $n\geq 1$, if the standard word $s_n$ in the language $St(\alpha_d)$ ends by a symbol $a$, then the word $a s_n^{d_{n+1}+2}$ is not a factor of $St(\alpha_d)$.
\end{lemma}

\section{The case of unique 1}\label{s:unique}
Let $r$ be a rotation word of the form $0^i10^{n-i}$. Clearly, if $0<i<n$, this word is generated by {\it all} pairs of Sturmian words of the form $(s10t,s01t)$ with $s10t,s01t\in S(n,\alpha)$ for some slope $\alpha$, where 
$|s|=i-1, |t|=n-i-1$. In particular, $s$ is a right special word in $St(i-1,\alpha)$, and $t$ is a left special word in $St(n-i-1,\alpha)$. Due to Claim \ref{cl:st_sp}, the shorter of the words $s$ and $t$ is determined by the longer one.

 Suppose first that $i\geq n-i$; then $s$ is longer than $t$ and determines all the word $s10t$. So, the number of pairs of Sturmian words giving $0^i10^{n-i}$ is equal to the number of (right) special Sturmian words of length $i-1$, that is, to 
$\sum_{p=1}^{i}\varphi(p)$ (see Claim \ref{cl:st_numbersp}). Symmetrically, if $i < n-i$, then this is $t$ that defines all the word $s10t$, and the number of such pairs is equal to the number of (left) special Sturmian words of length $n-i-1$, that is, to $\sum_{p=1}^{n-i}\varphi(p)$.

Now if $i=0$ and the rotation word is $10^n$, the pairs of Sturmian words generating it are all pairs of the form $0s$, $1s$. So, $s$ is a special Sturmian word of length $n-1$, and the number of such words is $\sum_{p=1}^{n}\varphi(p)$. Symmetrically, the number of pairs of Sturmian words generating the word $01^n$ is also the same. 

Summing up all the numbers above, we see that the $n+1$ rotation words of the form $0^i10^{n-i}$, where $0\leq i \leq n$, are generated by the following number of pairs of Sturmian words of length $n$:
\begin{equation}\label{e:f1}
f_1(n)=\begin{cases} 2 \sum_{i=k}^{2k}\sum_{p=1}^{i+1}\varphi(p), \mbox{~if~} n=2k+1, \\  
2 \sum_{i=k}^{2k-1}\sum_{p=1}^{i+1}\varphi(p) + \sum_{p=1}^{k}\varphi(p), \mbox{~if~} n=2k. \end{cases}
\end{equation}

Starting with $n=3$, the sums in $f_1(n)$ involve only special words of length at least 1. Exactly a half of them, namely, the left special words starting with 0 and symmetrically the right special words ending with 0, correspond to slopes less than 1/2. So, the $n+1$ rotation words of the form  $0^i10^{n-i}$ are generated by $f_1(n)/2$ pairs of Sturmian words.

Exactly the same total number of pairs (in fact, the pairs $(v,u)$, where $(u,v)$ are the pairs considered above) generate the $n+1$ rotation words of the form $1^i01^{n-i}$. Starting from $n=3$, it gives us exactly $f_1(n)$ pairs generating $2(n+1)$ rotation words.

\section{The case of several 1s}
Most of technical details of our result are hidden in the following
\begin{theorem}\label{ttt}
 Suppose a rotation word $w$ is generated by at least two different pairs  of Sturmian words of slope less than 1/2, and $w$ contains at least two 1s and at least two 0s. Then $w=0^i(10^l)^k10^j$ or $w=1^i(01^l)^k01^j$ for some $i,j\geq 0$, $l \geq 2$, $k \geq 1$, and the number of pairs generating $w$ is equal to $\varphi(l+1)/2$ if $i,j\leq l$ and $\varphi(l+1)$ otherwise.
\end{theorem}
This section is devoted to its proof which is based on the theory of standard Sturmian words and their construction with directive sequences (see Section \ref{s:prop}). 

\begin{lemma}\label{k1k2}
Suppose that a rotation word $r$ contains a factor $010^{k_1} 10^{k_2} 1$, where $k_2>k_1>0$, and is generated by ap pair $(u,v)$ of Sturmian words from some $St(\alpha_d)$. Then $k_1=l_n-1$ for some $n$ and $(k_2+1)\bmod l_n = l_{n-1}$.
\end{lemma}
{\sc Proof.} Clearly, the pair $(u,v)$ generating $r$ contains some factors $u'=10 u_1 10 u_2 1$ and 
$v'=01 u_1 01 u_2 0$ for some $u_1,u_22$ with $|u_1|=k_1-1$ and $|u_2|=k_2-1$. The words $u_1$ and $u_2$ are bispecial in $St(\alpha_d)$. So, $u_1=c_{n,d}$ for some $n>0$ and $0<d\leq d_n$. Note also that $u_2$ is left special in $St(\alpha_d)$ longer than $s_{n-1}$, so it starts from $s_{n-1}$ since there is only one left special word of each length in $St(\alpha_d)$.

Without loss of generality, suppose that $s_{n-2}$ ends by $10$; if it ends by $01$, in all the arguments below we should just consider $v'$ instead of $u'$.

Suppose that $d<d_n$; then $u'=10 s_{n-1}^d s_{n-2} s_{n-1} u''$; in particular, it means that the word $s_{n-1}^d s_{n-2} s_{n-1}$ can be extended to the left by $0$. On the other hand, since $d<d_n$, the same word can clearly be extended to the left by $s_{n-1}$, and thus by its last symbol $1$. We see that $s_{n-1}^d s_{n-2} s_{n-1}$ is left special; but it is not possible since $s_{n-2}s_{n-1}$ differs in two last symbols from $s_{n-1}s_{n-2}$, and thus  $s_{n-1}^d s_{n-2} s_{n-1}$ is not equal to the prefix $s_{n-1}^{d+1}s_{n-2}$ of $s_n$ of the same length. Since $s_{n-1}^{d+1}s_{n-2}$ is the only left special word of its length in $St(\alpha_d)$, we see that $0s_{n-1}^d s_{n-2} s_{n-1}$ and thus $u'$ are not in $St(\alpha_d)$, a contradiction. So, the case of $d<d_n$ is not possible, and thus $u_1=c_n$ for some $n$, and $k_1=l_n-1$.

Now recall that $u_2$ is longer than $u_1$; so, it is equal to $c_{N,D}$ for some $N>n$ and $0<D\leq d_N$. Suppose that $N>n+1$; then $u_2$ starts with $s_{n+1}$ and $u'= 10 s_n s_{n+1} u''$. As above, the word $s_n s_{n+1}$ can be extended to the left by $1$, which is the last symbol of $s_{n+1}$, and it is not left special since it is not equal to the only special word $s_{n+1}s_{n}$ of the same length, so, $0 s_n s_{n+1}$ and thus $u$ are not elements of $St(\alpha_d)$.

So, $N=n+1$, $u_2 01=s_n^d s_{n-1}$ for some $d\leq d_{n+1}$, and $|u_2|+2 \bmod l_n = l_{n-1}$, which was to be proved. \hfill $\Box$ 

\begin{lemma}\label{ddd}
 Consider two slopes $\alpha_d$, $\alpha_{d'}<1/2$ with corresponding directive sequences $(d_1,d_2,\ldots)$ and $(d_1',d_2',\ldots)$ and respective lengths $l_n$ and $l_n'$ of standard words. If $l_n=l_m'$ and $l_{n-1}=l_{m-1}'$ for some $m$ and $n$, then $n=m$ and $d_i=d_i'$ for all $i=1,\ldots n$. 
\end{lemma}
{\sc Proof.} By the construction, for all $i$ we have $l_{i-2}=l_{i} \bmod l_{i-1}$ and $d_i=\lfloor l_i:l_{i-1} \rfloor $. So, starting from $l_n=l_m'$ and $l_{n-1}=l_{m-1}'$, we can uniquely reconstruct $l_{n-2}=l_{m-2}'$, $l_{n-3}=l_{m-3}'$ etc. Note that here $d_1, d_1' >0$ since both slopes are less than 1/2. So, as soon as we get $l_{n-k}=l_{m-k}'=1$, we immediately see that $n=m=k$ and $(d_1,\ldots,d_n)=(d_1',\ldots,d_n')$. 
\hfill $\Box$

\begin{lemma}
 All rotation words arising from several pairs of Sturmian words of slope at most 1/2 and not containing two consecutive 1s are of the form $0^i(10^l)^k10^j$ for some $i,j\geq 0$, $k \geq 1$, $l\geq 2$.
\end{lemma}
{\sc Proof.} 
Let us prove that a word arising from two pairs of Sturmian words, of slopes $\alpha_d \neq \alpha_{d'}$, cannot contain a factor $w=10^{m_1}10^{m_2}1$ with $m_1\neq m_2$. Suppose it contains it. The proof is carried over for $m_2>m_1$; the opposite case can be proved by the argument that the set of Sturmian words, the set of rotation words and the procedure generating a rotation word from two Sturmian words are symmetric  under taking the mirror image.

First suppose that $1w$ is also a rotation word. Due to Lemma \ref{0011}, $1w$ appears from only one pair of Sturmian words, and this pair is $(0^{m_1+1}10^{m_2}1,010^{m_1}10^{m_2})$. These two words are Sturmian of the same slope only if $m_2=m_1+1$. But in this case, the word $0w=010^{m_1}10^{m_1+1}1$ due to Lemma \ref{k1k2} is also generated by only one pair, namely, by the pair $(10^{m_1}10^{m_1+1}1,010^{m_1}10^{m_2})$. So, $w$ arises only from the pair $(0^{m_1}10^{m_1+1}1,10^{m_1}10^{m_2})$, contradicting to our assumption.

So, if $w$ arises from several pairs of Sturmian words of slope at most 1/2, then so does $0w$. Since we suppose that $m_2>m_1$, we can apply Lemma \ref{k1k2}, according to which $m_1=l_n-1$ and $m_2\bmod m_1=l_{n-1}-1$, where $l_i$ are the lengths coming from the directive sequence for the language of the Sturmian words involved. Due to Lemma \ref{ddd}, the values of $l_n$ and $l_{n-1}$ uniquely determine the sequence $(d_1,d_2,\ldots, d_n)$, its length $n$ and thus the central words $u_1$ and $u_2$ such that $0w=0 u_1 10 u_2 1$ and $v=1 u_1 01 u_2 0$. So, we see that once again, $0w$ and thus $w$ arise from only one pair of Sturmian words of slope less than 1/2, a contradiction. 

We have proved that our word is of the form $0^i(10^l)^k10^j$ for some $i,j \geq 0$, $k \geq 1$ (since the case of a unique 1 is considered separately), and $l>0$. It remains to consider the case of $l=1$ and to see that the word $0^i(10)^k10^j$ is generated by only one pair of Sturmian words of slope less than 1/2 defined as follows: the central part $(10)^k1$ is given by the pair $u=0(10)^{k-1}1$, $v=1(01)^{k-1}0$, and the prefix and suffix zeros correspond to the common prefix $\cdots 1010$ and the common suffix $0101\cdots$ of the generating Sturmian words.
\hfill $\Box$ 

\begin{lemma}
 Each word $w=0^i(10^l)^k10^j$ with $i,j\leq l$, $l \geq 2$, $k \geq 1$ is generated by $\varphi(l+1)/2$ different pairs of Sturmian words of the same slope not greater than 1/2.
\end{lemma}
{\sc Proof.} Suppose for simplicity that $i,j>0$ and consider a pair $(u,v)$ such that $w=r(u,v)$. Clearly, $(u,v)=(s10c_1 10 \cdots, 10 c_k 10 p, s01c_1 01 \cdots, 01 c_k 01 p)$ for some central words $c_1,\ldots,c_k$ of length $l-1$ and some words $s,p$ with $|s|=i-1$, $|p|=j-1$. Since each Sturmian language  $St(\alpha)$ contains at most one central word of length $l-1$, we have $c_1=c_2=\cdots=c_k=c$. Moreover, the word $p$ is left special, and so it is a prefix of $c$, and the 
word $s$ is right special, so it is a suffix of $c$. So, the pair $(u,v)$ is uniquely determined by the central word $c$ and the parameter $i$. 

There exists $\varphi(l+1)$ central words of length $l-1$, and a half of them correspond to slopes less than 1/2. So, it remains to prove that for each central word $c$ and each power $k$ two words $u=s(10c)^k 10 p$ and $v=s(01c)^k p$, where $s$ is a suffix and $p$ is a prefix of $c$, appear in some Sturmian language of a given slope. 

Indeed, let $c10$ be equal to the standard word $s_n=s_{n-1}^{d_n}s_{n-2}$ in some language $St(\alpha_d)$. Then
$u$ is a factor of $s_n^{k+2}$. 
At the same time, $c01=s_n'=s_{n-1}^{d_n-1}s_{n-2} s_{n-1}$, so that $v$ is a factor of 
$(s')_n^{k+2}=s_{n-1}^{d_n-1}s_{n-2} s_n^{k+1}s_{n-1}$, which is in its turn a factor of $s_n^{k+3}$. So, 
taking $d_{n+1} \geq k+3$, we see that $u,v\in St(\alpha_d)$ for the directive sequence $d=(d_1,\ldots,d_n,k+3,\ldots)$, which was to be proved. 

If $i=0$, or $j=0$, or/and the standard word $s_n$ is equal to $c01$, not to $c10$, the proof is carried on similarly.
\hfill $\Box$

\begin{lemma}
 Each word $w=0^i(10^l)^k10^j$ with $i> l$ or $j>l$, $l \geq 2$, $k \geq 1$ is generated by $\varphi(l+1)$ different pairs of Sturmian words of the same slope not greater than 1/2.
\end{lemma}
{\sc Proof.} As above, if $w=r(u,v)$, then $u=s(10c)^k 10 p$ and $v=s(01c)^k 01 p$ for some central word $c$ of length $l-1$ and some words $s,p$ with $|s|=i-1$, $|p|=j-1$. We have a choice which of the words $c01$ and $c10$ is a standard word in the Sturmian language considered; suppose it is $c10=s_n$. Suppose also that $j>l$. Clearly, $p$ is left special and thus is a prefix of some standard word $s_N$, $N>n$. Suppose that $p$ is not a prefix of some power $s_n^d$ of $s_n$: it means that $p$ contains as a prefix the word $s_n^{d_{n+1}}s_{n-1}s_n$, or, more precisely, the word obtained from it by erasing the last symbol, since $s_{n-1}s_n$ differs from $s_ns_{n-1}$, which is a prefix of $s_n^2$, by the two last symbols (see Claim \ref{c:0110}). In particular, $p$ starts by $s_n^{d_{n+1}+1}$, and thus the suffix $0c10 p=0s_n p$ of $u$ starts with $0 s_{n}^{d_{n+1}+2}$, which is not an element of $St(\alpha_d)$ due to Lemma \ref{dn+2}, a contradiction. So, $p$ is a prefix of $s_n^d$ for some $d$. Symmetrically, if 
$i>l$, then $s$ is a suffix of the mirror image of $s_n^d$ for some $d$; by the way, this mirror image is equal to $01 (c01)^{d-1}c$. For $s_n=c01$, we should just consider $v$ instead of $u$ to prove the similar statements.

Note that since one of the parameters $i$ or $j$ is indeed greater than $l$, the cases of $s_n=c10$ and of  $s_n=c01$ are really different, which gives us $\varphi(l+1)$ cases: the total number of standard words of length $l+1$ is twice bigger than the number of central words of length $l-1$, that is, is equal to $2\varphi(l+1)$, but we are interested only in those of slope less than 1/2.

It remains to mention that for each standard word $s=cab$ of length $l+1$, where $a,b\in \{0,1\}$ $a \neq b$, and for all $d', d'', k \geq 0$, the  words $u'=(cba)^{d'}(cab)^{k+d''}$ and $v'=(cba)^{d'+k}(cab)^{d''}$, so that $w$ is a factor of $r(u',v')$, are factors of some language $St(\alpha_d)$. Indeed, let $s=s_n$ for the directive sequence $(d_1,\ldots,d_n)$; then $s=cab=s_{n-1}^{d_n}s_{n-2}$ and $cba=s_{n-1}^{d_n-1}s_{n-2}s_{n-1}$, so that $u'=(s_{n-1}^{d_n-1}s_{n-2}s_{n-1})^{d'}(s_{n-1}^{d_n}s_{n-2})^{k+d''}=s_{n-1}^{d_n-1}s_{n-2} s_n^{d'-1} s_{n-1} s_n^{k+d''}$, and $v'=$ $s_{n-1}^{d_n-1}s_{n-2} s_n^{d'+k-1} s_{n-1} s_n^{d''}$. If we take $d_{n+1}=D= \max\{d'+k,d''+k\}$, we see that both $u'$ and $v'$ are factors of $s_n^{D}s_{n-1}s_n^{D}$. So, both $u'$ and $v'$, and thus the pair of words based on the standard word $s$ of length $l+1$ and giving the rotation word $w$, are elements of the language $St(\alpha_d)$ for the directive sequence $d=(d_1,\ldots,d_n,D+2,\ldots)$, which completes the proof of the lemma. \hfill $\Box$

This lemma, in its turn, completes the proof of Theorem \ref{ttt}.

\section{Final computations}
To find the precise formula for $f(n+1)$ for $n \geq 3$, we should subtract from the bound $f_{pairs}(n)+2$, where $f_{pairs}(n)$ is found in Lemma \ref{l:pairs}, the number of pairs generating rotation words already obtained before. 

As it was shown in Section \ref{s:unique}, the $2(n+1)$ rotation words with only one symbol 1 or only one 0 are generated by $f_1(n)$ pairs (see \eqref{e:f1}). 

Now let us take into account the rotation words containing several 0s and several 1s. Their form is described in Theorem \ref{ttt}. Consider all rotation words of length $n+1$ of the form $0^i(10^l)^k10^j$, such that $l$ and $i$ are fixed and $k$ and $j$ are not. They are $\lfloor (n-i)/(l+1) \rfloor$; and taking all words with a given $l$ together, we see that they are
\[\sum_{i=0}^{n-(l+1)}\left \lfloor \frac{n-i}{l+1} \right \rfloor = \frac{1}{2}\left\lfloor \frac{n}{l+1} \right \rfloor (n-l+1+(n \bmod (l+1))).\]
In what follows, to make the formulas shorter, we will denote
\[
g(n,l)= n-l+1+(n \bmod (l+1)),
\]
so that the words of the form $0^i(10^l)^k10^j$ are $\frac{g(n,l) \lfloor n/(l+1)\rfloor}{2}$.

Each of these words is generated by $\varphi(l+1)$ pairs of Sturmian words, except for the $\min(l+1,n-l)$ words with $i,j\leq l$ which are generated by $\varphi(l+1)/2$ pairs each. So, for each $l \geq 2$ we should subtract from the sum the following function: 
{\small
\[
 f_2(n,l)=\left ( \frac{1}{2}\left \lfloor \frac{n}{l+1} \right\rfloor g(n,l)-\min(l+1,n-l)\right )(\varphi(l+1)-1)+\min(l+1,n-l)
\left ( \frac{\varphi(l+1)}{2}-1 \right ).
\]
}

The same function $f_2(n,l)$ corresponds to the words of the form $1^i(01^l)^k01^j$. So, to take into account all rotation words arising from several pairs and containing at least two 0s and at least two 1s, we should subtract from the upper bound the sum $2\sum_{l=2}^{n-1} f_2(n,l)$.

Summarizing the above arguments, we see that
\[f(n+1)=f_{pairs}(n)+2-f_1(n)+2(n+1)-2 \sum_{l=2}^{n-1} f_2(n,l).\]
This is exactly the statement of Theorem \ref{t:main} which was to be proved. \hfill $\Box$

\section{Acknowledgement}
 The authors are grateful to the participants of the working group in Saint-Di\'e des Vosges in April 2012 for stimulating discussions.

\end{document}